\documentclass[12pt]{amsart}
\usepackage{amsfonts, amsmath, amsthm}
\usepackage{epsfig}
\usepackage{psfrag}

\newtheorem{pro}{Proposition}[section]

\newtheorem{cnj}[pro]{Conjecture}

\newtheorem{quest}[pro]{Question}

\theoremstyle{definition}
\newtheorem{dfn}[pro]{Definition}

\theoremstyle{remark}

\title{Critical Heegaard Surfaces and Index 2 Minimal Surfaces} 
\date{August 1, 2001}
\address{Mathematics Department, University of Illinois at Chicago}
\email{bachman@math.uic.edu}
\author{David Bachman}
\begin{document}

\begin{abstract}
This paper contains the motivation for the study of {\it critical surfaces} in \cite{crit}. In that paper, the only justification given for the definition of this new class of surfaces is the strength of the results. However, when viewed as the topological analogue to index 2 minimal surfaces, critical surfaces become quite natural.
\end{abstract}
\maketitle

\section{Introduction.}

It is a standard exercise in 3-manifold topology to show that every manifold admits Heegaard splittings of arbitrarily high genus. Hence, a ``random" Heegaard splitting does not say much about the topology of the manifold in which it sits. To use Heegaard splittings to prove interesting theorems, one needs to make some kind of non-triviality assumption. The most obvious such assumption is that the splitting is minimal genus. However, this assumption alone is apparently very difficult to use. 

In \cite{cg:87}, Casson and Gordon define a new notion of triviality for a Heegaard splitting, called {\it weak reducibility}. A Heegaard splitting which is not weakly reducible, then, is said to be {\it strongly irreducible}. The assumption that a Heegaard splitting is strongly irreducible has proved to be much more useful than the assumption that it is minimal genus. In fact, in \cite{cg:87}, Casson and Gordon show that in a non-Haken 3-manifold, minimal genus Heegaard splittings {\it are} strongly irreducible. 

The moral here seems to be this: since the assumption of minimal genus is difficult to make use of, one should pass to a larger class of Heegaard splittings, which is still restrictive enough that one can prove non-trivial theorems. 

Now we switch gears a little. It is a Theorem of Riedemeister and Singer (see \cite{am:90}) that given two Heegaard splittings, one can always stabilize the higher genus one some number of times to obtain a stabilization of the lower genus one. However, this immediately implies that any two Heegaard splittings have a common stabilization of arbitrarily high genus. Hence, the assumption that one has a ``random" common stabilization cannot be terribly useful. What is of interest, of course, is the minimal genus common stabilization. As before though, the assumption of minimal genus has turned out to be very difficult to use. 

In this paper we will review the results of \cite{crit}, in which a new class of Heegaard splittings, called {\it critical}, was defined. The main result of that paper is that at least in the non-Haken case, this class includes the minimal genus common stabilizations. 

The term {\it critical} is defined via a 1-complex associated with any embedded, separating surface in a 3-manifold, which is reminiscent of the curve complex. The origin of the definition lies in the analogy between the study of the topology of embedded surfaces in 3-manifolds, and the study of minimal surfaces. This point of view was pioneered by mathematicians such as Hyam Rubinstein, Martin Scharlemann, and Abigail Thompson. We will make these analogies explicit here, and show how from this point of view, the study of critical surfaces is completely natural. 

For the sake of brevity, we will assume that the reader is familiar with the standard terminology of 3-manifold topology, that can be found in any introductory text. 

\section{Minimal surfaces}

Let $M$ be a 3-manifold, equipped with some Riemannian metric. Let $\Omega$ denote the space of all embedded surfaces in $M$, together with all surfaces that have been ``pinched" at finitely many points (at a point where a surface is ``pinched", there is a coordinate chart in which the surface looks locally like the graph of $z^2=x^2+y^2$), and finite collections of points. Hence, a path through the space $\Omega$ can be thought of as a continuous deformation of some surface, during which one might see some compressions (and ``de-compressions") happen. 

Let $A:\Omega \rightarrow \bf R$ denote the area function. What we are interested in is the critical points of $A$. Let $p$ be such a critical point, and let $D _p \nabla A$ denote the derivative of the gradient of $A$ at $p$. Let $\lambda _1, ..., \lambda _m$ denote the eigenvalues of $D _p \nabla A$ which are less than zero (we will only be interested in critical points in which the number of such eigenvalues is finite). Finally, let $n=\sum \limits _{i=1} ^m dim(V_i)$, where $V_i$ is the eigenspace corresponding to the eigenvector $\lambda _i$. If $S$ is the surface in $M$ which corresponds to the point $p$ of $\Omega$, then we say that $S$ is an {\it index $n$ minimal surface}. 

We now take a closer look at index 0, 1, and 2 minimal surfaces. 

\subsection{Index 0 Minimal Surfaces}
\label{ss:index0}

If $S$ is an index 0 minimal surface, and $p$ is the point of $\Omega$ corresponding to $S$, then it must be that $D _p \nabla A$ has only eigenvalues which are greater than or equal to zero. A much simpler way of saying this is that any perturbation of $S$ will increase its area, or keep it constant. In other words, $S$ is a local minimum for the area function, $A$. Hence, to locate an index 0 minimal surface, one simply needs to start with any surface, and ``flow downhill". That is, perturb it continuously in such a way so that its area decreases monotonically. Now, it may happen that what you end up with by doing this is a point, or a collection of points. Later we will encounter a topological restriction on the surfaces you can start with to avoid running into this problem. 

\subsection{Index 1 Minimal Surfaces}
\label{ss:index1}

To gain some intuition, it helps to visualize $\Omega$ as a 2-dimensional Euclidean space. In this case the graph of $A$ looks like a mountain range, and the index 1 critical points are the saddle points, i.e., the local maxima of the valleys. At such a point, one can travel through the valley forward or backwards, and ``go downhill", or one can leave the valley, and start to climb up the nearest mountain. 

Now, let's suppose we are faced with the task of finding an index 1 critical point. If we start with a random point and go downhill, then we will miss the index 1 points with probability 1. Instead, we can start with two index 0 points, and examine the paths which connect them. If we always keep the height of such a path as low as possible, then it is guaranteed to go through an index 1 point. 

The analogy is that you are a traveller in a mountain range. You want to get to your house, which is located in a pit, and you are starting out in some other pit. The catch is, you have some medical condition which makes you feel progressively more sick as your altitude increases. You would then choose to travel through the valleys, rather than climb over a mountain, even though that might be the shortest path. At some point in your journey, you will reach the highest point of some valley, and you will be at an index 1 critical point. 

One interesting note is that you may encounter several index 1 points along your way. If you always keep your altitude as low as possible, then in general your path will take you through a series of critical points which alternate between being index 0 and index 1. In general then, a manifold may contain a whole sequence of minimal surfaces, whose index alternates between 0 and 1.

\subsection{Index 2 Minimal Surfaces}
\label{ss:index2}

To find an index 1 minimal surface, we needed to start with two index 0 minimal surfaces, and connect them with an ``efficient" path. To find an index 2 minimal surface, we will now need to start with two efficient paths, and connect them with some kind of ``efficient" 1-parameter family {\it of paths}. Once again, by the term ``efficient" we mean that at all times we keep our area (i.e. our altitude in the mountain range) as small as possible. 

Let's look at our mountain range again for some intuition. It may be that you had two different choices of valleys to travel through when you were trying to get home in the previous subsection. Viewed from above, this would just look like two paths from one point to another. If we fill in the region between these two paths, then we are guaranteed to cover the peak of some mountain, an index 2 critical point.

\section{The Topological analogue of a minimal surface}

In this section, we look at analogues of index 0, 1, and 2 minimal surfaces, in the topological category. The main idea is to keep our space, $\Omega$, the same, but to change our area function to reflect only changes in topology. 

First, let $\Omega ^-$ denote subspace of $\Omega$ which consists of only the embedded surfaces of $M$. Now, let $A_T :\Omega ^- \rightarrow \bf Z$ denote the (continuous!) function defined by $A_T (S)=\sum \limits _{i=1} ^n (2-\chi (S_i))^2$, where $\{S_1,...,S_n\}$ are the components of $S$. Note that if $S$ is homeomorphic to $S^2$, then $A_T (S)=0$, reflecting the fact that in an irreducible 3-manifold (the kind we usually assume we are working in), every 2-sphere can be shrunk to a point. Also, if $S'$ is a surface obtained from the surface, $S$, by a compression, then $A_T (S')<A_T (S)$. 

\subsection{Incompressible Surfaces}

Recall that to find an index 0 minimal surface, you can try to start with some random embedded surface, and deform it through a path in $\Omega$ which monotonically decreases area. If you ever get ``stuck" at a surface with non-zero area, then you have found an index 0 minimal surface. 

To find an incompressible surface in a 3-manifold, $M$, one can start with some random surface, and try to compress it as much as possible. Each such compression decreases $A_T$, and so in some sense we are doing the exact analogue of what we did in subsection \ref{ss:index0}. If the process ever terminates in anything other than a union of 2-spheres (i.e. in anything for which the function, $A_T$, is non-zero), then one has found an incompressible surface. Hence, for many purposes it is useful to think of incompressible surfaces as the topological analogue of index 0 minimal surfaces. 

\subsection{Strongly Irreducible Surfaces}

To look for some kind of toplogical surface that is the appropriate analogue of an index 1 minimal surface, we follow the strategy of subsection \ref{ss:index1}. That is, we start with two incompressible surfaces (the analogues of index 0 minimal surfaces), and look at the paths from one to the other in $\Omega$ for which the function, $A_T$, is always as small as possible. This is precisely the strategy of Scharlemann and Thompson from \cite{st:94}, in which they show that every irreducible 3-manifold contains a strongly irreducible Heegaard splitting for some submanifold cobounded by (possibly empty) incompressible surfaces (in the case of empty incompressible surfaces, their techniques just produce a strongly irreducible Heegaard splitting for the entire manifold). Hence, we will view strongly irreducible Heegaard splittings as the appropriate analogues of index 1 minimal surfaces. 

The analogy is really quite good. In subsection \ref{ss:index1}, we saw that an efficient path connecting two index 0 minimal surfaces may contain not just one index 1 minimal surface, but a whole sequence of minimal surfaces whose index alternates between 0 and 1. In the topological category, the aforementioned work of Scharlemann and Thompson shows that every irreducible 3-manifold admits an alternating sequence of incompressible surfaces, and strongly irreducible Heegaard splittings.

\subsection{Critical surfaces}

We now come to the main point of this paper, which was to motivate the study of a new class of topological surfaces. These new surfaces arise naturally as the appropriate analogues of index 2 minimal surfaces. Since index $n$ minimal surfaces correspond to critical points, we have chosen the name {\it critical} for our new class. The precise definition of a critical surface will be given in the next section. First, we say in what sense they are analogous to index 2 minimal surfaces. 

Recall from subsection \ref{ss:index2} that index 2 minimal surfaces arise when perturbing one ``efficient" path through $\Omega$ to another. Along the way, we will see a sequence of paths that can be described like this: odd elements of the sequence go through an alternating sequence of index 0 and index 1 minimal surfaces. Even elements are similar, except that in place of exactly one index 1 minimal surface we see an index 2 minimal surface. 

In the topological category, the strategy is exactly the same. We look at two different sequences of surfaces which alternate between being incompressible and strongly irreducible. We then try to ``connect" these sequences with intermediate sequences in such a way so that the function, $A_T$, is as small as possible at all times. This is precisely the strategy of \cite{crit}. In that paper, we show that the sequences alternate as follows: odd sequences contain an alternating sequence of incompressible and strongly irreducible surfaces. Even sequences are similar, except that in place of exactly one strongly irreducible surface there is a critical surface.

\section{The definition of a Critical Surface.}

To facilitate the definition of a critical surface, we first define a 1-complex for each embedded, orientable, closed, separating surface in a 3-manifold, $M$. Suppose $F$ is such a surface. If $D$ and $D'$ are compressing disks for $F$, then we say $D$ is equivalent to $D'$ if there is an isotopy of $M$ taking $F$ to $F$, and $D$ to $D'$ (we do allow $D$ and $D'$ to be on opposite sides of $F$). 

We now define a 1-complex, $\Gamma (F)$. For each equivalence class of compressing disk for $F$, there is a vertex of $\Gamma (F)$. Two (not necessarily distinct) vertices are connected by an edge if there are representatives of the corresponding equivalence classes on opposite sides of $F$, which intersect in at most a point. A vertex of $\Gamma (F)$ is said to be {\it isolated} if it is not the endpoint of any edge.

For example, if $F$ is the genus 1 Heegaard splitting of $S^3$, then there is an isotopy of $S^3$ which takes $F$ back to itself, but switches the sides of $F$. Such an isotopy takes a compressing disk on one side of $F$ to a compressing disk on the other. Hence, $\Gamma (F)$ has a single vertex. However, there are representatives of the equivalence class that corresponds to this vertex which are on opposite sides of $F$, and intersect in a point. Hence, there is an edge of $\Gamma (F)$ which connects the vertex to itself.

\begin{dfn}
If we remove the isolated vertices from $\Gamma (F)$ and are left with a disconnected 1-complex, then we say $F$ is {\it critical}.
\end{dfn}

Equivalently, $F$ is critical if there exist two edges of $\Gamma(F)$ that can not be connected by a 1-chain.

\section{Results about critical surfaces}

We now state some of the main results about critical Heegaard surfaces from \cite{crit}. The first few Lemmas of that paper build up to the following Theorem:

\medskip
\noindent {\bf Theorem 4.6.} {\it Suppose $M$ is an irreducible 3-manifold with no closed incompressible surfaces, and at most one Heegaard splitting (up to isotopy) of each genus. Then $M$ does not contain a critical Heegaard surface.} 
\medskip

The remainder of \cite{crit} is concerned with the converse of this Theorem. That is, we answer precisely when a (non-Haken) 3-manifold {\it does} contain a critical Heegaard surface. 

The main technical theorem which starts us off in this direction is:

\medskip
\noindent {\bf Theorem 5.1.} {\it Let $M$ be a 3-manifold with critical surface, $F$, and incompressible surface, $S$. Then there is an incompressible surface, $S'$, homeomorphic to $S$, such that every loop of $F \cap S'$ is essential on both surfaces. Furthermore, if $M$ is irreducible, then there is such an $S'$ which is isotopic to $S$.}
\medskip

Note that this Theorem was already known to be true if one were to replace the word ``critical" with either ``incompressible" or ``strongly irreducible". This is just more evidence for the naturality of critical surfaces.

As immediate corollaries to this, we obtain:

\noindent {\bf Corollary 5.7.} {\it A reducible 3-manifold does not admit a critical Heegaard splitting.}

\medskip
\noindent {\bf Corollary 5.8.} {\it Suppose $M$ is a 3-manifold which admits a critical Heegaard splitting, such that $\partial M \ne \emptyset$. Then $\partial M$ is essential in $M$.}
\medskip

It is this last corollary which we combine with a considerable amount of new machinery (all motivated by the analogy with index 2 minimal surfaces) to yield:

\medskip
\noindent {\bf Theorem 6.1.} {\it Suppose $F$ and $F'$ are distinct strongly irreducible Heegaard splittings of some closed 3-manifold, $M$. If the minimal genus common stabilization of $F$ and $F'$ is not critical, then $M$ contains an incompressible surface.}
\medskip

We actually prove a slightly stronger version of this Theorem, that holds for manifolds with non-empty boundary. 

Compare Theorem 6.1. to that of Casson and Gordon \cite{cg:87}: If the minimal genus Heegaard splitting of a 3-manifold, $M$, is not strongly irreducible, then $M$ contains an incompressible surface. Once again, the parallels between these two Theorems is yet more evidence for the naturality of critical surfaces.

\section{A Conjecture}

The relationships between index 0 and 1 minimal surfaces and incompressible and strongly irreducible Heegaard splittings seem to be much deeper than mere analogy. For instance, in \cite{fhs:83}, Freedman, Hass and Scott show that any incompressible surface can be isotoped to be a least area surface. Such surfaces are index 0 minimal surfaces. In \cite{pr:87}, Pits and Rubinstein show that strongly irreducible surfaces can always be isotoped to index 1 minimal surfaces. This motivates us to make the following conjecture:

\begin{cnj}
Any critical surface can be isotoped to be an index 2 minimal surface. 
\end{cnj}

In \cite{crit2}, we prove a Piecewise-Linear analogue of this.

\section{A metric on the space of strongly irreducible Heegaard splittings}

We now show how our results lead to a natural metric on the space of strongly irreducible Heegaard splittings of a non-Haken 3-manifold. The author believes that it would be of interest to understand this space better. 

First, given a critical surface, $F$, we can define a larger 1-complex, $\Lambda(F)$, that contains $\Gamma(F)$ as follows: the vertices of $\Lambda(F)$ are equivalence classes of loops on $F$, where two loops are considered equivalent if there is an isotopy of $M$ taking $F$ to $F$, and one loop to the other. There is an edge connecting two vertices if there are representatives of the corresponding equivalence classes which intersect in at most a point. Recall that a vertex of $\Gamma(F)$ corresponds to an equivalence class of compressing disks for $F$. Thus, we can identify each vertex of $\Gamma(F)$ with the vertex of $\Lambda(F)$ which corresponds to the boundary of any representative disk. 

Now, suppose $e_1$ and $e_2$ are two edges in $\Gamma (F)$. Define $d(e_1,e_2)$ to be the minimal length of any chain connecting $e_1$ to $e_2$ in $\Lambda(F)$. Now, given two components, $C_1$ and $C_2$, of $\Gamma (F)$, we can define $d(C_1,C_2)$, the {\it distance} between $C_1$ and $C_2$, as $\min \{d(e_1,e_2)|e_1$ is an edge in $C_1$, and $e_2$ is an edge of $C_2\}$.

Finally, suppose $H_1$ and $H_2$ are strongly irreducible Heegaard splittings of a 3-manifold, $M$, and $F$ is their minimal genus common stabilization. As $F$ is a stabilization of $H_i$, it can be isotoped so that between $F$ and $H_i$ there is a compression body, $W_i$, and so that there are compressing disks for $F$, $D_i \subset W_i$, and $E_i \subset cl(M -W_i)$, such that $|D_i \cap E_i|=1$. Each pair, $(D_i,E_i)$ corresponds to some edge of $\Gamma (F)$. In \cite{crit}, we show that the edge corresponding to $(D_1,E_1)$ is in a component, $C_1$, of $\Gamma (F)$ which is different than the component, $C_2$, containing the edge corresponding to $(D_2,E_2)$, and that $C_1$ and $C_2$ were independent of our exact choices of $D_i$ and $E_i$. We can therefore define the {\it distance} between $H$ and $H'$ as $d(C_1,C_2)$. 

\begin{quest}
Can the distance between strongly irreducible Heegaard splittings be arbitrarily high? If not, is there a bound in terms of the genera of the splittings, or perhaps a universal bound?
\end{quest}

\begin{quest}
Is there an algorithm to compute the distance between two given strongly irreducible Heegaard splittings?
\end{quest}

\begin{quest}
Is there a relationship between the distance between two strongly irreducible Heegaard splittings, and the number of times one needs to stabilize the higher genus one to obtain a stabilization of the lower genus one?
\end{quest}

\begin{quest}
Is there a relationship between the distance between two strongly irreducible Heegaard splittings, and the distances of each individual splitting, in the sense of Hempel \cite{hempel:01}?
\end{quest}

\bibliographystyle{plain}
\bibliography{kyoto}

\end{document}